\documentclass[12pt]{amsart}
\usepackage{amssymb,amsmath}

\calclayout

\newtheorem{theorem}{Theorem}
\newtheorem{corollary}[theorem]{Corollary}
\newtheorem{proposition}[theorem]{Proposition}
\newtheorem{lemma}[theorem]{Lemma}

\newcommand{\cC}{{\mathcal C}}
\newcommand{\cH}{{\mathcal H}}
\newcommand{\cI}{{\mathcal I}}
\newcommand{\cO}{{\mathcal O}}
\newcommand{\cU}{{\mathcal U}}

\newcommand{\mP}{{\mathbb P}}

\newcommand{\oB}{\bar{B}}
\newcommand{\oG}{\bar{G}}
\newcommand{\oP}{\bar{P}}
\newcommand{\oS}{\bar{S}}
\newcommand{\oT}{\bar{T}}

\newcommand{\fg}{{\mathfrak g}}
\newcommand{\fm}{{\mathfrak m}}

\newcommand{\diag}{{\rm diag}}
\newcommand{\Chow}{{\rm Chow}}
\newcommand{\Grass}{{\rm Grass}}
\newcommand{\Hilb}{{\rm Hilb}}
\newcommand{\Aut}{{\rm Aut}}
\newcommand{\End}{{\rm End}}
\newcommand{\Hom}{{\rm Hom}}
\newcommand{\Spec}{{\rm Spec}}
\newcommand{\Stab}{{\rm Stab}}
\newcommand{\supp}{{\rm supp}}

\title{Group completions via Hilbert schemes}
\author{Michel~Brion}
\address{Universit\'e de Grenoble I\\
D\'epartement de Math\'ematiques\\
Institut Fourier, UMR 5582 du CNRS\\
38402 Saint-Martin d'H\`eres Cedex, France}
\email{Michel.Brion@ujf-grenoble.fr}

\date{}
\begin{document}

\begin{abstract}
Let $X$ be a projective variety, homogeneous under a linear algebraic
group. We show that the diagonal of $X$ belongs to a unique irreducible
component $\cH_X$ of the Hilbert scheme of $X\times X$. Moreover,
$\cH_X$ is isomorphic to the ``wonderful completion'' of the connected
automorphism group of $X$; in particular, $\cH_X$ is non-singular. We
describe explicitly the degenerations of the diagonal in $X\times X$,
that is, the points of $\cH_X$; these subschemes of $X\times X$ are
reduced and Cohen-Macaulay. 
\end{abstract}

\maketitle

\section*{Introduction}

Let $G$ be an adjoint semisimple algebraic group over an algebraically
closed field. Then $G$ admits a canonical smooth completion
$\oG$, such that: (a) the action of $G\times G$ by left and right
multiplication on $G$ extends to $\oG$, (b) the boundary $\oG - G$ is
a union of smooth irreducible divisors intersecting transversally
along the unique closed $G\times G$-orbit, and (c) the partial
intersections of these boundary divisors are exactly the orbit
closures. This ``wonderful'' completion was constructed by De Concini
and Procesi \cite{DP} in characteristic zero, and then by Strickland
\cite{St} in arbitrary characteristics, via representation theory.

\medskip

In this article, we obtain algebro-geometric realizations of $\oG$, as
follows. Choose a parabolic subgroup $P$ of $G$ and consider the
projective variety $G/P=X$. Regard the diagonal in $X\times X$ as a
point of the Hilbert scheme $\Hilb(X\times X)$. The group $G\times G$
acts on $\Hilb(X\times X)$ via its natural action on $X\times X$; if
$G$ acts faithfully on $X$, then the $G\times G$-orbit of the diagonal
is isomorphic to $G$. We show that the closure of this
orbit is isomorphic to $\oG$ (Theorem \ref{main}). If moreover $G$
is the full connected automorphism group of $X$, this realizes $\oG$
as the irreducible component $\cH_X$ of $\Hilb(X\times X)$ through the
diagonal (Lemma \ref{component}).  

\medskip

In the case where $G$ is the projective linear group ${\rm PGL}(n+1)$,
the wonderful completion $\oG$ is the classical ``space of complete 
collineations'', and the most natural choice for $X$ is the projective
space $\mP^n$. In fact, the realization of $\oG$ as the orbit
closure of the diagonal in the Hilbert scheme of $\mP^n\times\mP^n$ is
due to M.~Thaddeus \cite{T}, in the more general setting of ``complete
collineations of a given rank''. The approach of Thaddeus is inspired
by Kapranov's work on the moduli space of stable punctured curves of
genus $0$, see \cite{K1}, \cite{K2}. It proceeds through several other
interesting constructions of $\oG$, either as an iterated blow-up or
as a quotient, that do not seem to extend to other groups.

\medskip

For this reason, we follow an alternative approach based on
results of \cite{B} and \cite{BP}. A flat family over $\oG$
was constructed there, whose general fibers are the 
$G\times G$-translates of the diagonal in $X\times X$. The resulting
morphism $\oG\to\Hilb(X\times X)$ turns out to be an isomorphism over
its image (Theorem \ref{main}). As a consequence, we describe the
fibers over $\oG$ of the universal family of $\Hilb(X\times X)$, that
is, the scheme-theoretic degenerations of the diagonal. These
subschemes turn out to be reduced and Cohen-Macaulay. They are
obtained by moving a union of Schubert varieties in $X\times X$ by the
diagonal action of a Levi subgroup, see Proposition \ref{fibers}. The
case where $X$ is the full flag variety of $G$ turns out to be nicer,
e.g., all degenerations are Gorenstein. (This does not even extend to
$X=\mP^2$, see the remark at the end of Section 3.)

\medskip

This approach also yields algebro-geometric realizations of the
``wonderful symmetric varieties'' of \cite{DP}, \cite{DS}; moreover,
the Hilbert scheme may be replaced everywhere with the Chow variety,
in characteristic zero (Theorem \ref{chow} and Corollary
\ref{chowtwist}). This generalizes results of Thaddeus \cite{T}
concerning  complete quadrics and complete skew forms.

\medskip

As another generalization, one may consider a compact K\"ahler
manifold $X$ and its connected automorphism group $G$; then the
$G\times G$-orbit closure $\cH_X$ of the diagonal in the Douady space
of $X\times X$ is again an equivariant compactification of $G$, which
may be of interest (this construction is used in \cite{Sn} to
compactify the $G$-action on $X$). The description of $\cH_X$ is easy
if $X$ is homogeneous: for, by a theorem of Borel and Remmert (see
\cite{A} 3.9), $X$ is then the product of its Albanese variety $A(X)$,
a complex torus, with a flag variety $Y$ as above. And one checks that
$\cH_{Y\times A(X)}\cong\cH_Y\times\cH_{A(X)}\cong\cH_Y\times A(X)$. 

\medskip

More generally, for any integer $m\geq 2$, the irreducible component
of the ``small'' diagonal in the Hilbert or Douady space of $X^m$ is
an equivariant compactification of the homogeneous space $G^m/G$,
where $G=\Aut^0(X)$ is embedded diagonally in $G^m$.  
Are some of these completions smooth ? In the case where 
$G={\rm PGL}(n+1)$, could they be related to those constructed by
Lafforgue \cite{L} ? 

\medskip

This work is organized as follows. Section 1 contains preliminary
constructions and remarks concerning connected automorphism groups
of arbitrary projective varieties. Our main result, Theorem
\ref{main}, is stated and proved in Section 2, and then extended to
wonderful symmetric varieties. As an application, all
degenerations of the diagonal in $G/P\times G/P$ are described in
Section 3. In Section 4, we adapt our constructions and results to
Chow varieties, assuming that the characteristic of $k$ is
zero. We obtain a slightly stronger version of Theorem \ref{main}
under this assumption, with a much simpler proof.

\medskip

We use \cite{Sp} as a general reference for linear algebraic groups,
and \cite{Ko} for Hilbert schemes.

\bigskip

\noindent
{\sl Acknowledgements}. Many thanks to to Etienne Ghys and Laurent
Manivel for useful discussions, and special thanks to Bill Haboush for
his inspiring questions during a conference on algebraic groups held
at Aarhus in June 2000.

\section{Completions of connected automorphism groups}

In this section, we consider an arbitrary projective algebraic
variety $X$ over an algebraically closed field $k$. Recall that the
automorphism group $\Aut(X)$ is an algebraic group scheme over $k$
with Lie algebra $\Hom(\Omega^1_X,\cO_X)$, 
where $\Omega^1_X$ denotes the sheaf of K\"ahler differentials of $X$
(see \cite{Ko} I.2.16.4). Let $G$ be a closed connected
subgroup scheme of $\Aut(X)$, with Lie algebra $\fg$. We shall
assume that $G$ is smooth; this assumption is always fulfilled in
characteristic zero.
  
The group $G\times G$ acts on $X\times X$, and thus on its Hilbert
scheme $\Hilb(X\times X)$. Let $\cH_{X,G}$ be the $G\times G$-orbit
closure in $\Hilb(X\times X)$ of the diagonal, $\diag(X)$; we endow
$\cH_{X,G}$ with its reduced subscheme structure. 

\begin{lemma}\label{completion}
With preceding notation, $\cH_{X,G}$ is a projective 
$G\times G$-equivariant completion of the homogeneous space 
$(G\times G)/\diag(G)$.
\end{lemma}

\begin{proof}
Clearly, the isotropy group of $\diag(X)$ in $G\times G$ equals
$\diag(G)$, as a set. Thus, it suffices to check that the isotropy Lie
algebra of $\diag(X)$ is $\diag(\fg)$. To see this, recall that the
Zariski tangent space to $\Hilb(X\times X)$ at $\diag(X)$ is
$$
T_{\diag(X)}\Hilb(X\times X)= 
\Hom(\cI_{\diag(X)}/\cI^2_{\diag(X)},\cO_{\diag(X)})
=\Hom(\Omega_X^1,\cO_X).
$$
Moreover, the differential of the morphism 
$$
G\times G\to\Hilb(X\times X),~(g,h)\mapsto (g,h)\cdot\diag(X)
$$
at the identity element of $G\times G$, is the map
$$
\fg\times\fg\to\Hom(\Omega_X^1,\cO_X),~(x,y)\mapsto x-y.
$$
\end{proof}

Note that $(G\times G)/\diag(G)$ is isomorphic to $G$ via 
$(g,h)\mapsto gh^{-1}$; moreover, the subscheme $(g,h)\cdot\diag(X)$
is the graph of $hg^{-1}=(gh^{-1})^{-1}$. Thus, the closed points of
$\cH_{X,G}$ are the graphs of elements of $G$, together with their
limits as closed subschemes.

If $\Aut(X)$ is smooth, then we can take $G=\Aut^0(X)$;  
we then set $\cH_{X,G}=\cH_X$. By the preceding argument, the Zariski
tangent spaces at $\diag(X)$ of $\cH_X$ and of $\Hilb(X\times X)$
coincide; we thus obtain

\begin{lemma}\label{component}
With preceding notation, $\diag(X)$ is a non-singular point of
$\Hilb(X\times X)$. Moreover, $\cH_X$ is the irreducible component 
of $\Hilb(X\times X)$ through $\diag(X)$. 
\end{lemma}

\medskip

The $G\times G$-action on $\cH_{X,G}$ can be interpreted in terms of
two well-known operations on $\Hilb(X\times X)$. The first one is the
``convolution'' product $*$, defined by 
$$
Z_1 * Z_2 =p_{13*}(p_{12}^*Z_1\cap p_{23}^*Z_2)
$$
where $p_{ij}:X\times X\times X\to X\times X$ denotes the projection
to the $(i,j)$-factor, and $\cap$ denotes scheme-theoretic
intersection. The second one is the ``transposition'', the involution
of $\Hilb(X\times X)$ induced by the involution $(x,y)\mapsto(y,x)$ of
$X\times X$; this involution will be denoted $Z\mapsto Z^{-1}$. 

For every automorphism $g$ of $X$, with graph $\Gamma_g$, we have
$\Gamma_g^{-1}=\Gamma_{g^{-1}}$. Moreover, left (resp.~right)
convolution with $\Gamma_g$ is the action of $(g^{-1},id)$
(resp.~$(id,g)$) on $\Hilb(X\times X)$.  Thus, the transposition
extends to $\cH_{X,G}$ the inverse map on $G$, and the $G\times G$-action
on $\cH_{X,G}$ can be viewed as left and right convolution with graphs of
elements of $G$. 

However, $\cH_{X,G}$ is generally not invariant under
convolution (thus, convolution is not a morphism.) For example, if
$X=\mP^1$ and $G={\rm PGL}(2)$, then $\cH_{X,G}=\cH_X$ contains the
reduced subscheme $Z_x=\mP^1\times\{x\}\cup\{x\}\times\mP^1$ for every
$x\in\mP^1$, and $Z_x * Z_x$ equals $\mP^1\times\mP^1$.

\medskip

Finally, we shall need a ``twisted'' version of $\cH_{X,G}$. Let
$\sigma$ be an automorphism of the algebraic group $G$, such that the
fixed point subgroup scheme 
$$
G^{\sigma}=\{g\in G~\vert~\sigma(g)=g\}
$$
is smooth (e.g., $G$ is linear and $\sigma$ is semisimple). Consider
another variety $Y$ endowed with an action of $G$ and with an
isomorphism
$$
f:X\to Y,
$$ 
such that $f(g\cdot x)=\sigma(g)\cdot f(x)$ for all 
$(g,x)\in G\times X$. Let 
$$
\Gamma_f\subseteq X\times Y
$$ 
be the graph of $f$. The stabilizer of $\Gamma_f$ under
the diagonal action of $G$ on $X\times Y$ is $G^{\sigma}$. Let
$\cH_{f,G}$ be the $G$-orbit closure of $\Gamma_f$ in the Hilbert
scheme of $X\times Y$, then $\cH_{f,G}$ is a projective
$G$-equivariant completion of the homogeneous space $G/G^{\sigma}$;
its points are graphs of $G$-conjugates of $f$, together with their
limits as subschemes.

Note that the isomorphism 
$$
id\times f^{-1}:X\times Y\to X\times X
$$
maps $\Gamma_f$ to $\diag(X)$. Thus, $\cH_{f,G}$ embeds into
$\cH_{X,G}$ as a closed subvariety, invariant under the ``twisted''
action of $G$ defined by
$$
g * Z=(g,\sigma^{-1}(g))\cdot Z.
$$
The restriction to the open orbit 
$G\cdot\Gamma_f \cong G/G^{\sigma}$ 
of both embeddings identifies to the map
$$
G/G^{\sigma}\to G,~gG^{\sigma}\mapsto g\sigma^{-1}(g^{-1}),
$$
equivariant for the $G$-action on $G/G^{\sigma}$ by left
multiplication, and for the twisted $G$-action on $G$.

\medskip

\noindent
{\sl Example 1.} Consider a finite dimensional $k$-vector space
$V$ endowed with a nondegenerate quadratic form $q$; we assume that
the characteristic of $k$ is not $2$. This defines a
symmetric isomorphism $V\to V^*$ and hence an isomorphism
$f:\mP(V)\to\mP(V^*)$, together with an
involutive automorphism $\sigma$ of ${\rm PGL}(V)$ mapping every
element to the inverse of its adjoint with respect to $q$. Clearly,
we have $f(g\cdot x)=\sigma(g)\cdot f(x)$ for all 
$(g,x)\in{\rm PGL}(V)\times\mP(V)$, and $\sigma$ is semisimple with
fixed point subgroup the projective orthogonal group, 
${\rm PO}(V,q)$. Thus, we obtain an equivariant completion
$\cH$ of the homogeneous space ${\rm PGL}(V)/{\rm PO}(V,q)$. The
latter is the space of non-singular quadrics in $\mP(V)$, whereas
$\cH$ is the closure of the locus of graphs of symmetric isomorphisms
$\mP(V)\to\mP(V^*)$, in the Hilbert scheme of
$\mP(V)\times\mP(V^*)$. 

It was shown by Thaddeus that $\cH$ is isomorphic to the ``space of
complete quadrics'' in $\mP(V)$, see \cite{T} \S 10. Moreover,
replacing the quadratic form $q$ by a non-degenerate skew form yields
the ``space of complete skew forms'', see \cite{T} \S 11. Both spaces
are classical examples of wonderful symmetric varieties, see
\cite{DP}; in fact, Thaddeus' result will be generalized to all
wonderful symmetric varieties in the next section.

\section{Wonderful completions via Hilbert schemes of flag varieties}

In this section, we consider a connected linear algebraic group $G$
and a parabolic subgroup $P$ such that $G$ acts faithfully on
$G/P=X$. Then $G$ is semisimple adjoint. Moreover, the group scheme
$\Aut^0(X)$ is smooth and is semisimple adjoint as well; this
group equals $G$ as a rule, and all exceptions are known (see
\cite{De}). Thus, we may consider the equivariant completion
$\cH_{X,G}$ (resp.~$\cH_X$); we aim at proving that it is the
wonderful completion of $G$ (resp.~of $\Aut^0(X)$). For this, we
review results of \cite{B} and \cite{BP}. 

Let $\oP$ be the closure of $P$ in the wonderful completion
$\oG$. Since $\oP$ is invariant under the action of $P\times P$ on
$\oG$, we can form the associated fiber bundle
$$
p: G\times G\times^{P\times P}\oP\to 
(G\times G)/(P\times P)=X\times X,
$$
a locally trivial $G\times G$-equivariant fibration with fiber
$\oP$. Moreover, the inclusion of $\oP$ into
$\oG$ extends uniquely to a $G\times G$-equivariant map
$$
\pi:G\times G\times^{P\times P}\oP\to\oG.
$$
The product map
$$
p\times\pi:G\times G\times^{P\times P}\oP \to X \times X \times \oG
$$
is a closed immersion; its image is the ``incidence variety''
$$
\{(gP,hP,x)~\vert~ (g,h)\in G\times G,~x\in(g,h)\oP\}.
$$
Thus, the fibers of $\pi$ identify to closed subschemes of $X\times X$
via $p_*$; this identifies the fiber at the identity element of $G$,
to $\diag(X)$. By \cite{B} 1.6, $\pi$ is equidimensional; moreover, by  
\cite{BP} \S 7, $\oP$ is Cohen-Macaulay, and the fibers of $\pi$ are
reduced. It follows that $\pi$ is flat, with Cohen-Macaulay fibers. 

By the universal property of the Hilbert scheme and the
definition of $\cH_{X,G}$, we thus obtain a morphism 
$$
\varphi:\oG\to \cH_{X,G}
$$
that maps every $\gamma\in\oG$ to the subscheme $p_*(\pi^*\gamma)$ of
$X\times X$. 

\begin{theorem}\label{main}
With preceding notation, $\varphi$ is an isomorphism. As a
consequence, the restriction to $\cH_{X,G}$ of the universal family
over $\Hilb(X\times X)$ identifies to $\pi$, so that its fibers are
reduced Cohen-Macaulay subschemes of $X\times X$.
\end{theorem}

\begin{proof}
Choose opposite Borel subgroups $B$, $B^-$ of $G$, with common
torus $T$ and unipotent radicals $U$, $U^-$. By the Bruhat
decomposition, the product map 
$$
U\times T\times U^-\to G
$$ 
is an open immersion. By \cite{DP} and \cite{St}, this map extends to
an open immersion
$$
U\times \oT_0\times U^- \to \oG_0,
$$ 
where $\oT_0$ is an affine open subset of the closure of $T$ in $\oG$. 
Moreover, $\oT_0$ is isomorphic to an affine space, and it meets
every $G\times G$-orbit in $\oG$ along a unique $T\times T$-orbit. As a
consequence, $\oG_0$ is an $B\times B^-$-invariant open affine subset
of $\oG$, and $\oT_0$ meets the closed $G\times G$-orbit (isomorphic to
$G/B^-\times G/B$) at the unique point $z_0=B^-/B^-\times B/B$. 

It suffices to show that $\varphi:\oG\to\cH_{X,G}$ is bijective, and
restricts to an isomorphism over $\oG_0$. For the latter assertion 
(the main point of the proof), we shall show that every subscheme
$Z\in\varphi(\oG_0)$ intersects transversally certain Schubert
varieties in $X\times X$. Mapping $Z$ to its common points with these
Schubert varieties, and to the Zariski tangent spaces of $Z$ at these
points, will yield appropriate morphisms from $\varphi(\oG_0)$ to $U$,
$U^-$ and $\oT_0$, and hence an inverse to $\varphi\vert_{\oG_0}$.

We need more notation, and results from \cite{BP}. Choose for $B$ a
Borel subgroup of $P$; let $L$ be the Levi subgroup of $P$ containing
the maximal torus $T$. Let $W$ be the Weyl group of $(G,T)$ and let
$\Phi$ be the root system of $(G,T)$, with the subset $\Phi^+$ of
positive roots defined by $B$ and the corresponding basis $\Delta$.
Every $\alpha\in\Delta$ defines a simple reflection $s_{\alpha}\in
W$. The $s_{\alpha}$ ($\alpha\in\Delta)$ generate $W$; this defines
the length function $\ell$ on that group. Let $w_{\Delta}$ be the
unique element of maximal length; then 
$B^- = w_{\Delta} B w_{\Delta}$. Let $I$ be the subset of $\Delta$
consisting of all simple roots of $(L,T)$, and let $\Phi_I$ be the
corresponding root subsystem of $\Phi$. The Weyl group of $(L,T)$ is
the Weyl group $W_I$ of the root system $\Phi_I$; it is generated by
the $s_{\alpha}$, $\alpha\in I$. 
Let 
$$
W\!^I=\{w\in W~\vert~w(I)\subseteq\Phi^+\},
$$
then $W\!^I$ consists of all elements of $W$ of minimal length in
their right $W_I$-coset; it is a system of representatives of the
quotient $W/W_I$. The latter identifies to the double coset space
$B\backslash G/P$, via $w\mapsto BwP/P$. 

Recall that the $T$-fixed points in $G/P$ are the $e_w=wP/P$ for 
$w\in W\!^I$, and the $B$-orbits (the Bruhat cells) are the
$$
C_w=B\cdot e_w\cong U/U\cap {^w\!P}
$$ 
(where ${^w\!P}=wPw^{-1}$); the Schubert varieties are their closures 
$$
X_w=\overline{B\cdot e_w}.
$$ 
We shall also need the opposite Bruhat cells 
$$
C^-_w=B^-\cdot e_w,
$$ 
and the opposite Schubert varieties 
$$
X^-_w=\overline{B^-\cdot e_w}=w_{\Delta}X_{w_{\Delta}ww_I}
$$
(note that $W\!^I$ is invariant under $w\mapsto w_{\Delta}ww_I$).

By \cite{B} 1.6, we have
$$
p(\pi^{-1}(z_0))=\bigcup_{w\in W\!^I} X^-_w\times X_w
$$
as sets. Thus, the same holds for the scheme $\varphi(z_0)=Z_0$, 
since it is reduced.

\begin{lemma}\label{bijective}
The isotropy group of $Z_0$ in $G\times G$ equals 
$B^-\times B$ (as sets). Moreover, $\varphi$ is bijective.
\end{lemma}

\begin{proof}

The isotropy group $\Stab_{G\times G}(Z_0)$ contains $B^-\times B$; 
thus, it can be written as $Q^-\times Q$ for parabolic subgroups $Q^-$ 
containing $B^-$, and $Q$ containing $B$. Since $Q$ is connected, it
stabilizes all Schubert varieties $X_w$ ($w\in W\!^I$). As a
consequence, $Q$ stabilizes all Bruhat cells.

If $Q\neq B$, then we can choose $\alpha\in\Delta$ such that
$s_{\alpha}$ has a representative in $Q$. Let $w\in W\!^I$, then
$Bs_{\alpha}BwP=BwP$, so that $s_{\alpha}w\in wW_I$ by the Bruhat
decomposition. Thus, $w^{-1}(\alpha)\in\Phi_I$. Since this 
holds for all $w$ in a system of representatives of $W/W_I$, it
follows that the orbit $W\alpha$ is contained in $\Phi_I$. Therefore,
$\Phi_I$ contains the intersection of $\Phi$ with the linear span of
$W\alpha$. As a consequence, $P$ contains a non-trivial closed normal
subgroup of $G$. But this contradicts the assumption that $G$ acts
faithfully on $G/P$. 

Thus, $Q=B$; likewise, $Q^-=B^-$. In other words, the restriction of
$\varphi$ to the closed $G\times G$-orbit is bijective. Using
\cite{Kn} \S 7, it follows that $\varphi$ is the normalization of its
image; in particular, $\varphi$ is finite. Moreover,
by \cite{DP} and \cite{DS}, the $G\times G$-isotropy group of every
point of $\oG$ is connected, and the conjugacy class of that group
determines the orbit uniquely (an explicit description of these
isotropy groups will be recalled in Section 3). Therefore, $\varphi$
is bijective.
\end{proof}

As a consequence, 
$\varphi(\oG_0)=(U\times U^-)\cdot \varphi(\oT_0)$ 
is an open subset of $\varphi(\oG)=\cH_{X,G}$. We now investigate the
structure of this subset, in a succession of lemmas.

\begin{lemma}\label{transversal}
For all $Z\in\varphi(\oT_0)$ and $w\in W\!^I$, the scheme-theoretic
intersection of $Z$ with $X_w\times X^-_w$ is supported at the 
unique point $(e_w,e_w)$, and this intersection is transversal.
\end{lemma}

\begin{proof}
Note that $\varphi(1)=\diag(X)$ contains the $T\times T$-fixed point
$(e_w,e_w)$; therefore, this point belongs to every $Z\in\varphi(\oT)$.

We first determine the intersection 
$Z_0\cap(X_w\times X^-_w)$. It is invariant under
$T\times T$, and contains $(e_w,e_w)$. Let $(e_u,e_v)$ be another
fixed point. Then, by the description of $Z_0$, there exists $x\in
W\!^I$ such that: $e_u\in X^-_x\cap X_w$ and $e_v\in X_x\cap X^-_w$. 
Thus, $u\leq ww_I$ for the Bruhat ordering on $W$, so that $u\leq w$. 
Moreover, $e_u\in w_{\Delta}X_{w_{\Delta}xw_I}$, so that $x\leq u$. 
Likewise, we have $w\leq v\leq x$. Therefore, we must have $x=w=u=v$, 
and $(e_w,e_w)$ is the unique $T\times T$-fixed point of the set
$Z_0\cap(X_w\times X^-_w)$. It follows that $(e_w,e_w)$ is
the unique point of that set. Since $X_w$ and $X^-_w$ intersect
transversally at $e_w$, we see that $Z_0$ and 
$X_w\times X^-_w$ intersect transversally at $(e_w,e_w)$.

To extend this to all $Z\in\varphi(\oT_0)$, consider the universal
family $\cU$ of $\Hilb(X\times X)$, and its pullback $\cU_{\oT_0}$ to
$\oT_0$. This is a closed subscheme of $X\times X\times\oT_0$; the
scheme-theoretic intersection 
$$
\cU_{\oT_0}\cap (X_w\times X^-_w\times\oT_0)
$$ 
is invariant under the natural action of $T\times T$, and the second
projection
$$
q:\cU_{\oT_0}\cap (X_w\times X^-_w\times\oT_0)\to\oT_0
$$
is equivariant and proper. By the preceding discussion, the
scheme-theoretic fiber of $q$ at is the closed point
$(e_w,e_w,z_0)$. Since each $T\times T$-orbit closure in $\oT_0$
contains $z_0$, it follows that all fibers of $q$ are finite. Thus,
$q$ is finite, and $\cU_{\oT_0}\cap (X_w\times X^-_w\times\oT_0)$ is
affine. Now a version of Nakayama's lemma implies that $q$ is an
isomorphism.
\end{proof}

By Lemma \ref{transversal} and the structure of $\oG_0$, every
$Z\in\varphi(\oG_0)$ intersects $X_w\times X^-_w$ transversally, at
a unique point of $C_w\times C^-_w$; let $p_w(Z)$ be that point.

\begin{lemma}\label{points}
Every $p_w:\varphi(\oG_0)\to C_w\times C^-_w$ is a 
$B\times B^-$-equivariant morphism, and the image of the product map
$$
\prod_{w\in W\!^I} p_w:\varphi(\oG_0)\to
\prod_{w\in W\!^I} C_w\times C^-_w
$$
is closed and isomorphic to $U\times U^-$. 
\end{lemma}

\begin{proof}
Consider again the universal family of $\Hilb(X\times X)$, and its
restriction $\cU_{\varphi(\oG_0)}$ to $\varphi(\oG_0)$. Then the
scheme-theoretic intersection  
$$
\cU_{\varphi(\oG_0)}\cap (X_w\times X^-_w\times\varphi(\oG_0))
$$ 
projects isomorphically to $\varphi(\oG_0)$, for every fiber is a
point. It follows that $p_w$ is a morphism; it is clearly
$B\times B^-$-equivariant.  

Thus, the product morphism $\prod_{w\in W\!^I} p_w$ is equivariant as
well. To show that its image is isomorphic to $U\times U^-$, it
suffices to check that the map
$$
\iota:U\to \prod_{w\in W\!^I} C_w,~u\mapsto (u\cdot e_w)_{w\in W\!^I}
$$
is a closed immersion. For this, identify each Bruhat cell $C_w$
with the homogeneous space $U/U\cap {^w\!P}$. Arguing as in the proof of
Lemma \ref{bijective}, we see that the intersection of the subgroups
$U\cap {^w\!P}$ ($w\in W\!^I$) is trivial, and that the same holds for the 
intersection of their Lie algebras. Thus, the orbit map $\iota$ is an
immersion. Moreover, its image is closed, as an orbit of an unipotent
group acting on an affine variety.
\end{proof}

By Lemma \ref{transversal}, every $Z\in\varphi(\oT_0)$ contains
$(e_w,e_w)$ as a non-singular point, and the Zariski tangent space 
$$
T_{(e_w,e_w)} Z=t_w(Z)
$$
is transversal to 
$$
T_{(e_w,e_w)} (X_w\times X^-_w)=T_{(e_w,e_w)} (C_w\times C^-_w).
$$ 
Let $\Grass_w$ be the Grassmanian of linear subspaces
of $T_{(e_w,e_w)} (X\times X)$, and let $\Grass_{w,0}$ be the open
affine subset consisting of those subspaces that are transversal to
$T_{(e_w,e_w)} (C_w\times C^-_w)$.

\begin{lemma}\label{spaces}
Every $t_w:\varphi(\oT_0)\to\Grass_{w,0}$ is a morphism, and the image
of the product map
$$
\prod_{w\in W\!^I} t_w:\varphi(\oT_0)\to 
\prod_{w\in W\!^I}\Grass_{w,0}
$$
is isomorphic to $\oT_0$.
\end{lemma}

\begin{proof}
Consider once more the universal family of $\Hilb(X\times X)$, and its
restriction $\cU_{\varphi(\oT_0)}$ to $\varphi(\oT_0)$. Let $\fm$ be
the ideal sheaf of the closed point $(e_w,e_w)$ of $X\times X$ and 
consider the scheme-theoretic intersection
$$
\cU_{\varphi(\oT_0)}\cap 
(\Spec(\cO_{X\times X}/\fm^2)\times \varphi(\oT_0)),
$$
with projection map $q$ to $\varphi(\oT_0)$. Then the fiber of $q$ at
$Z$ equals $\Spec(\cO_Z/\fm^2\cO_Z)$; its length is finite and
independent of $Z$, since $(e_w,e_w)$ is a non-singular point of that
scheme. Thus, $q$ is finite and flat; it follows that the assignment
$$
Z\mapsto (\fm\cO_Z/\fm^2\cO_Z)^*=t_w(Z)
$$ 
is a morphism. Note that $t_w$ maps $\diag(X)$ to the diagonal of
$T_{e_w}X$.

For the second assertion, fix $w\in W\!^I$ and set 
$$
V=T_{e_w}X,~ V_+=T_{e_w}C_w \text{ and } V_-=T_{e_w}C^-_w.
$$ 
Then $V$ is the direct sum of its subspaces $V_+$ and
$V_-$. Thus, every subspace of $V\times V=T_{(e_w,e_w)} X\times X$,
transversal to the subspace 
$V_+\times V_-=T_{(e_w,e_w)}C_w\times C_w^-$, can be written as a graph
$$
\{(v_- + f_+(v_+ + v_-), v_+ + f_-(v_+ + v_-)~\vert~
v_+\in V_+,~v_-\in V_- \},
$$
with uniquely defined linear maps $f_{\pm}:V\to V_{\pm}$. This
defines isomorphisms
$$
\Grass_{w,0}\cong\Hom(V,V_+)\times\Hom(V,V^-)\cong\End(V),
$$
mapping $\diag(V)$ to the identity. 
Moreover, for every $t\in T$, the point 
$$\displaylines{
(t,1)\cdot\diag(V) =\{(t\cdot (v_- + v_+),v_+ + v_-)
~\vert~ v_+\in V_+,~v_-\in V_- \}
\hfill\cr\hfill
=\{(v_- + t \cdot v_+, v_+ + t^{-1}\cdot v_-)
~\vert~ v_+\in V_+,~v_-\in V_- \}
\cr}$$
is mapped to the endomorphism 
$$
\diag(t,t^{-1}):V\to V,~v_+ + v_-\mapsto t\cdot v_+ + t^{-1}\cdot v_-.
$$
Since each weight space of $V$ is one-dimensional, it follows that
the eigenvalues of $\diag(t,t^{-1})$ are non-zero regular functions
on $\varphi(\oT_0)$, eigenvectors of $T$ (acting on the left) with
weights: the $T$-weights of $V_-$ and the opposites of the $T$-weights
of $V_+$. Let $\Phi_w$ be the set of weights of these functions; for
$\chi\in\Phi_w$, let $f_{\chi}$ be the corresponding function. Then
$f_{\chi}$ is an eigenvector of $T$ (acting on the right) of weight
$-\chi$.

The set of weights of $V=T_{e_w}G/P$ (resp.~$V_+$, $V_-$) equals
$w(\Phi^- - \Phi_I)$ (resp. $\Phi^+\cap w(\Phi^- - \Phi_I)$; 
$\Phi^-\cap w(\Phi^- - \Phi_I)$). Thus, we have
$$
\Phi_w = \Phi^-\cap w(\Phi-\Phi_I) = \Phi^- - w(\Phi_I).
$$
Moreover, for every simple root $\alpha$, there exists $w\in W\!^I$
such that $-\alpha$ belongs to $\Phi^- - w(\Phi_I)$, by
the proof of Lemma \ref{bijective}; and the functions $f_{\alpha}$
($\alpha\in\Delta$) generate the coordinate ring of $\oT_0$, by
\cite{DS} \S 3. Thus, the composition of $\varphi$ with 
$\prod_{w\in W\!^I} t_w$ is a closed immersion.
\end{proof}

Lemmas \ref{points} and \ref{spaces} imply that the restriction
$$
\varphi:\oG_0\cong U\times U^-\times\oT_0\to\varphi(\oG_0)
$$ 
is an isomorphism. Since $\varphi$ is bijective and $\oG_0$ meets all
$G\times G$-orbits in $\oG$, it follows that $\varphi$ is an
isomorphism.
\end{proof}

\medskip

Next we extend Theorem \ref{main} to symmetric spaces; for this, we
assume that the characteristic of $k$ is not $2$. Let $\sigma$ be
an automorphism of order $2$ of $G$, with fixed point subgroup
$G^{\sigma}$. By \cite{DP} and \cite{DS}, the adjoint symmetric space
$G/G^{\sigma}$ admits a canonical $G$-equivariant completion
$\overline{G/G^{\sigma}}$, a ``wonderful symmetric variety''. 
We shall need the following realization of $\overline{G/G^{\sigma}}$
in $\oG$, obtained by Littelmann and Procesi in characteristic zero
(see \cite{LP} 3.2).

\begin{lemma}\label{fixed}
The map $G/G^{\sigma}\to G$, $gG^{\sigma}\mapsto \sigma(g)g^{-1}$
extends to a closed embedding of $\overline{G/G^{\sigma}}$ into $\oG$.
\end{lemma}

\begin{proof}
Since the approach of \cite{LP} 3.2 does not extend to arbitrary
characteristics in a straightforward way, we provide an alternative
argument. First we review some results from \cite{DS}.
  
Let $P$ be a parabolic subgroup of $G$ such that $^{\sigma}\!P$ and $P$
are opposite, and that $P$ is minimal for this property. Choose a
maximal $\sigma$-split subtorus $S$ of $P$, a maximal torus $T$ of $P$
containing $S$, and a Borel subgroup $B$ of $P$ containing $T$. Then
$T$ is $\sigma$-stable, so that $\sigma$ acts on the character group
of $T$, and on the root system $\Phi$; the opposite Borel subgroup
$B^-$ is contained in $^\sigma\!P$, and contains $R_u({^\sigma\!P})$.
Moreover, the natural map 
$R_u(P)\times S/S^{\sigma}\to G/G^{\sigma}$
extends to an open immersion 
$$
R_u(P)\times (\overline{S/S^{\sigma}})_0\to \overline{G/G^{\sigma}},
$$
where $(\overline{S/S^{\sigma}})_0$ is isomorphic to affine space
where $S$ acts linearly with weights $\alpha-\sigma(\alpha)$,
$\alpha\in\Delta$.  

The isomorphism $U\times U^-\times\oT_0\to \oG_0$ restricts to a
closed immersion
$$
\iota: R_u(P)\times \oS_0\to\oG_0,
~(g,\gamma)\mapsto (g,\sigma(g))\cdot\gamma,
$$
where $\oS_0=\oS\cap\oT_0$. Note that $\oS_0$ is isomorphic to
$(\overline{S/S^{\sigma}})_0$, equivariantly for the action of $S$ on
$\oS_0$ by $(g,\gamma)\mapsto g^2\cdot\gamma$, and for the natural
action of $S$ on $(\overline{S/S^{\sigma}})_0$. Moreover, 
$\iota(R_u(P)\times S)$ is contained in the image of $G/G^{\sigma}$ in
$G$ (for every $g\in S$ can be written as
$\gamma^2=\gamma\sigma(\gamma)^{-1}$ for some $\gamma\in S$). 

Let $\oG_{\sigma}$ be the closure in $\oG$ of the image of
$G/G^{\sigma}$, and let $\oG_{\sigma,0}=\oG_{\sigma}\cap\oG_0$. Then
$\iota$ induces an isomorphism
$$
R_u(P)\times\oS_0\to\oG_{\sigma,0}.
$$
It follows that the rational map 
$\overline{G/G^{\sigma}}\to \oG_{\sigma}$ is defined on
$(\overline{G/G^{\sigma}})_0$ and maps it isomorphically to
$\oG_{\sigma,0}$. Since $(\overline{G/G^{\sigma}})_0$ meets all
$G$-orbits in $\overline{G/G^{\sigma}}$, this rational map is an
isomorphism.
\end{proof}

We return to the situation where $P$ is a parabolic subgroup of $G$
such that $G$ acts faithfully on $G/P$. Then $\sigma(P)$ is another
parabolic subgroup of $G$, and $\sigma$ induces an isomorphism 
$$
f:G/P\to G/\sigma(P),
$$ 
such that $f(g\cdot x)=\sigma(g)\cdot f(x)$ for all 
$(g,x)\in G\times G/P$. Thus, we can consider the equivariant completion
$\cH_{f,G}$ of the symmetric space $G/G^{\sigma}$, constructed at the
end of Section 1. The discussion in that section, together with
Theorem \ref{main}, yields immediately the following result.

\begin{corollary}\label{twist}
With preceding notation, $\cH_{f,G}$ is isomorphic to the wonderful
completion of $G/G^{\sigma}$.
\end{corollary}

\section{The degenerations of the diagonal of a flag variety}

We still consider an adjoint semisimple group $G$ and a parabolic
subgroup $P$, such that $G$ acts faithfully on $G/P=X$. In this
section, we shall describe the points of $\cH_{X,G}$ viewed as closed
subschemes of $X\times X$, that is, the partial degenerations of the
diagonal (the total degenerations being the points of the closed
orbit). By Theorem \ref{main}, this amounts to describing the
set-theoretical fibers of the map
$$
\pi:G\times G\times^{P\times P}\oP\to\oG,
$$
embedded into $G/P\times G/P$ via the projection
$$
p:G\times G\times^{P\times P}\oP\to G/P\times G/P.
$$
Moreover, by equivariance, it suffices to determine these fibers
$p(\pi^{-1}(x))$ at representatives $x$ of the (finitely many) orbits
of $G\times G$ in $\oG$. This will be achieved in Proposition
\ref{fibers} below, in terms of combinatorics of Weyl groups;
geometric applications will be given after the proof of that
Proposition. 

We use the notation $B$, $B^-$, $T$, $\Phi$, $\Phi^+$, $\Delta$, $W$,
$\ell$, $w_{\Delta}$, $W_I$, $W\!^I$, $X_w$, $X^-_w$ introduced in
the proof of Theorem \ref{main}. For every (possibly empty) subset $J$
of $\Delta$, let $\lambda_J$ be the unique one-parameter subgroup of
$T$ such that
$$
\langle\lambda_J,\alpha\rangle=
\begin{cases}0&\text{ if }\alpha\in J;\\
1&\text{ if }\alpha\in\Delta \setminus J.\\
\end{cases}
$$
Let $x_J$ be the limit in $\oG$ of $\lambda_J(t)$ as 
$t\to 0$. By \cite{DS} \S 3 (see also the Appendix in
\cite{B}), the $x_J$ ($J\subseteq\Delta$) are a system of
representatives of the $G\times G$-orbits in $\oG$; note that
$\lambda_{\Delta}=0$, so that $x_{\Delta}$ is the identity element of
$G$. This sets up an order-preserving bijection between subsets $J$ of
$\Delta$ and $G\times G$-orbits $\cO_J$ (ordered by inclusion of their 
closures). In particular, the closed orbit is $\cO_{\emptyset}$.

Every $\lambda_J$ determines two opposite parabolic subgroups $P_J$,
$P_J^-$ of $G$, where $P_J$ (resp.~$P_J^-$) consists of those $g\in G$
such that $\lambda_J(t)g\lambda_J(t)^{-1}$ has a limit in $G$ as $t\to 0$
(resp.~$t\to\infty$). Note that $P_J$ contains $B=P_{\emptyset}$,
whereas $P_J^-$ contains $B^-=Q_{\emptyset}$. The common Levi subgroup
$L_J=P_J\cap P_J^-$ is the centralizer of the image of $\lambda_J$; its
root system is $\Phi_J$. Moreover, $B_J=B\cap L_J$ is a Borel subgroup
of $L_J$. 

Now the isotropy group scheme $\Stab_{G\times G}(x_J)$ is smooth; it
is the semi-direct product of the unipotent radical of 
$P_J^-\times P_J$ with $\diag(L_J)(C_J\times C_J)$, where $C_J$ is the
center of $L_J$. In particular, the isotropy group of $x_{\emptyset}$
is $B^-\times B$.

Since $C_J=\{t\in T~\vert~\alpha(t)=1~\forall\alpha\in J\}$ and $J$ is
part of the basis $\Delta$ of the character group of $T$, the group
$C_J$ is connected; as a consequence, $\Stab_{G\times G}(x_J)$ is
connected as well. 

Let $I$ be the subset of $\Delta$ such that $P=P_I$. Let
$$
^J\!W\!^I=\{w\in W~\vert~w(I)\subseteq\Phi^+  \text{ and }
w^{-1}(J)\subset\Phi^+\}.
$$
By \cite{Bo} IV.1, Exercice 3, this subset of $W$ consists of all
elements $w$ of minimal length in their double coset $W_JwW_I$;
moreover, it is a system of representatives of the double coset space 
$W_J\backslash W/W_I$, or, equivalently, of 
$P_J\backslash G/P=P_J\backslash X$ by \cite{Bo} IV.2.6, Proposition 2.

For every $w$ in $^J\!W\!^I$, we have $L_J\cap {^w\!B}=B_J$; thus,
$L_J\cap {^w\!P}$ is a parabolic subgroup of $L_J$.

\begin{proposition}\label{fibers}
With preceding notation, every irreducible component of the fiber
$p(\pi^{-1}(x_J))$ can be written as
$$
Z_w=\diag(L_J)\cdot(w_J X^-_{w_Jw}\times X_w)
$$
for a unique $w$ in $^J\!W\!^I$. Moreover,
$w_J X^-_{w_Jw}\times X_w$ is invariant under the diagonal action
of $L_J\cap {^w\!P}$, and the variety
$$
\tilde Z_w = L_J\times^{L_J\cap {^w\!P}}
(w_J X^-_{w_Jw}\times X_w)
$$
maps birationally to $Z_w$ under the natural morphism
$\rho_w:\tilde Z_w\to Z_w$.
\end{proposition}

\begin{proof}
As a first step, we show the equality in $\oG$:
$$
\oP\cap\overline{\cO_J}=\bigcup_{w\in ^I\!W\!^J}
\overline{(P\times P)(w,w)\cdot x_J}
$$
(decomposition into irreducible components). For this, note that we
have
$$
\oP=\overline{Bw_IB}=\overline{Bw_Iw_{\Delta}B^-}w_{\Delta}
=(1,w_{\Delta})\cdot\overline{Bw_Iw_{\Delta}B^-}.
$$
Moreover, by \cite{B} Theorem 2.1, the irreducible components of
$\overline{Bw_Iw_{\Delta}B^-}\cap\overline{\cO_J}$ are the closures
$$
\overline{(B\times B^-)(w_Iw_{\Delta}v,v)\cdot x_J},
$$
where $v\in W$ satisfies
$$
v\in W\!^J \text{ and }
\ell(w_Iw_{\Delta})=\ell(w_Iw_{\Delta}v)+\ell(v).
$$
Thus, we obtain the following decomposition into irreducible
components: 
$$
\oP\cap\overline{\cO_J}=\bigcup
\overline{(B\times B)(w_Iw_{\Delta}v,w_{\Delta}v)\cdot x_J},
$$
the union over all $v$ as above. Since every irreducible component of 
$\oP\cap\overline{\cO_J}$ is invariant under $P\times P$, we
can rewrite this as
$$
\oP\cap\overline{\cO_J}=\bigcup
\overline{(P\times P)(w_{\Delta}v,w_{\Delta}v)\cdot x_J}.
$$
Set $w=w_{\Delta}v$. Then, since $W\!^J$ is invariant under the map
$v\mapsto w_{\Delta}vw_J$, we have
$$
v\in W\!^J\Leftrightarrow ww_J\in W\!^J.
$$
This in turn amounts to: $w$ is the unique element of maximal length
in its right $W_J$-coset. On the other hand, we have
$$
\ell(w_Iw_{\Delta})=\ell(w_Iw_{\Delta}v)+\ell(v)
\Leftrightarrow \ell(w_Iw)=\ell(w)-\ell(w_I)
\Leftrightarrow w_Iw\in ^I\!W,
$$
which is equivalent to: $w$ is the element of maximal length in its
left $W_I$-coset. So we have proved that
$$
\oP\cap\overline{\cO_J}=\bigcup
\overline{(P\times P)(w,w)\cdot x_J},
$$
the union over all $w\in W$ of maximal length in their right
$W_J$-coset and in their left $W_I$-coset. But every irreducible
component $\overline{(P\times P)(w,w)\cdot x_J}$ depends only on the
double coset $W_IwW_J$ (for the isotropy group of $x_J$ contains
$\diag(L_J)$). Moreover, this double coset contains a unique element
of maximal length (for $W_IwW_J$ is the set of $B\times B$-orbits in
$Pw P_J$, and the latter contains a unique open 
$B\times B$-orbit). Thus, the preceding union is over the set
$W_I\backslash W/W_J$, or, equivalently, over $^I\!W\!^J$. This
completes the first step of the proof.

As a second step, we show that
$$
p(\pi^{-1}(x_J))=\bigcup_{w\in ^J\!W\!^I} 
\overline{\Stab_{G\times G}(x_J)\cdot (e_w,e_w)}
$$
(decomposition into irreducible components). For this, given
$(g,h)\in G$, note that 
$$(gP,hP)\in p(\pi^{-1}(x_J))\Leftrightarrow
(g^{-1},h^{-1})\cdot x_J\in\oP\cap (G\times G)\cdot x_J.
$$
Note also that $^I\!W\!^J$ and $^J\!W\!^I$ are exchanged by 
$w\mapsto w^{-1}$. The assertion follows from these remarks, together
with the first step. 

Now recall that
$$
\Stab_{G\times G}(x_J)=\diag(L_J)(R_u(P_J^-)C_J\times R_u(P_J)C_J).
$$
Moreover, for every $w$ in $^J\!W\!^I$, we have
$$
R_u(P_J)C_JwP=R_u(P_J)wP=BwP,
$$
for $B=R_u(P_J)B_J$ and $B_J\subseteq {^w\!P}$. As a consequence, 
$BwP$ is invariant under left multiplication by the group 
$L_J\cap {^w\!P}$. Likewise, since $w_JB^-w_J=R_u(P_J^-)B_J$, we have
$$
R_u(P_J^-)C_JwP=R_u(P_J^-)wP=w_JB^-w_JwP,
$$
and this subset is $L_J\cap {^w\!P}$-invariant. Thus,
$$
\Stab_{G\times G}(x_J)\cdot (e_w,e_w)=
\diag(L_J)\cdot(w_J C^-_{w_Jw}\times C_w).
$$
Together with the second step, this proves all assertions of the
Proposition, except for birationality of $\rho_w$. But $\tilde Z_w$
contains 
$$\displaylines{
\tilde Z_w^0=L_J\times^{L_J\cap {^w\!P}}
(R_u(P_J^-)\cdot e_w \times R_u(P_J)\cdot e_w
\hfill\cr\hfill
=L_J\times^{L_J\cap {^w\!P}}
((R_u(P_J^-)/R_u(P_J^-)\cap {^w\!P})\times (R_u(P_J)/R_u(P_J)\cap {^w\!P}))
\cr}$$
as an open subset, mapped under $\rho_w$ onto 
$$
Z_w^0=\Stab_{G\times G}(x_J)\cdot(e_w,e_w).
$$ 
We shall show that the restriction 
$$
\rho_w:\tilde Z_w^0\to Z_w^0
$$ 
is an isomorphism.
 
For this, we describe the (set-theoretical) isotropy group of
$(e_w,e_w)$ in $\Stab_{G\times G}(x_J)$, that is, 
$\Stab_{{^w\!P}\times{^w\!P}}(x_J)$. Since 
$$
\Stab_{G\times G}(x_J)=
(R_u(P_J^-)\times R_u(P_J))\cdot (C_J\times C_J)\cdot \diag(L_J),
$$
this isotropy group is contained in 
$$
({^w\!P}\cap P_J^-)\times ({^w\!P}\cap P_J) =
({^w\!P}\cap R_u(P_J^-))({^w\!P}\cap L_J)\times
({^w\!P}\cap R_u(P_J))({^w\!P}\cap L_J),
$$
and contains 
$({^w\!P}\cap R_u(P_J^-))\times({^w\!P}\cap R_u(P_J))$. 
It follows that
$$
\Stab_{{^w\!P}\times {^w\!P}}(x_J)=
(R_u(P_J^-)\cap {^w\!P})\times(R_u(P_J)\cap {^w\!P})
\cdot (C_J\times C_J) \cdot \diag(L_J\cap {^w\!P}).
$$
The isotropy Lie algebra of $(e_w,e_w)$ in 
$\Stab_{\fg\times \fg}(x_J)$ is described similarly; it follows that 
$$
Z^0_w\cong \Stab_{G\times G}(x_J)/\Stab_{{^w\!P}\times {^w\!P}}(x_J)
$$
and that $\rho_w:\tilde Z_w^0\to Z_w^0$ is bijective and separable,
hence an isomorphism.
\end{proof}

\begin{corollary}\label{reducible}
With preceding notation, the set of irreducible components of
$\pi^{-1}(x_J)$ is in bijection with the set of $P_J$-orbits 
in $X$. As a consequence, the fibers of $\pi$ at all boundary points
of $\bar{G}$ are reducible. 
\end{corollary}

\begin{proof}
The first assertion is a direct consequence of Proposition
\ref{reducible}. If $\pi^{-1}(x_J)$ is irreducible, then $P_J$ acts
transitively on $X$, and hence its unipotent radical acts trivially. 
Since $G$ acts faithfully on $X$, it follows that $P_J=G$, that is,
$J=\emptyset$.
\end{proof}

\medskip

\noindent
{\sl Example 2.} Let $X$ be the full flag variety of $G$, that is,
$P=B$, or, equivalently, $I$ is empty. Then $^J\!W\!^I$ equals
$^J\!W$, and $L_J\cap {^w\!P}=B_J$, so that the statement of
Proposition \ref{fibers} simplifies slightly as follows:
$$
p(\pi^{-1}(x_J))= \bigcup_{w\in ^J\!W} 
\diag(L_J)\cdot (w_J X^-_{w_Jw}\times X_w)
$$
(decomposition into irreducible components), and the map
$$
\rho:L_J\times^{B_J} \bigcup_{w\in ^J\!W}
w_J X^-_{w_Jw}\times X_w\to p(\pi^{-1}(x_J))
$$
is birational. 

\medskip

\noindent
{\sl Example 3.} Let $X=\mP^n=\mP(k^{n+1})$, then 
$G={\rm PGL}(n+1)$. Let $(v_1,\ldots,v_{n+1})$ be the standard basis
of $k^{n+1}$; let $T$ (resp.~$B$, $B^-$) be the image in $G$ of the
group of diagonal (resp.~upper triangular, lower triangular) matrices
in ${\rm GL}(n+1)$. The Weyl group $W$ is the permutation group of
$1,\ldots,n+1$, and the simple reflections are the transpositions
$(i,i+1)$ where $1\leq i\leq n$. Identifying the simple roots with the
corresponding reflections, we have
$$
I=\{(2,3),\ldots,(n,n+1)\} \text{ and }
W^I=\{w\in W~\vert~ w(2)<w(3)< \cdots <w(n+1)\}.
$$
Thus, every $w\in W^I$ is uniquely determined by $w(1)$; for 
$1\leq i\leq n+1$, we denote $w_i$ the element of $W^I$ such that
$w_i(1)=i$. Then $e_{w_i}$ is the line $kv_i$, and the
Schubert variety $X_{w_i}$ is simply the subspace
$\mP(k v_1+\cdots+k v_i)$, whereas the opposite Schubert variety 
$X^-_{w_i}$ equals $\mP(k v_i+\cdots+k v_{n+1})$.

Let 
$$
J=\{(j_1,j_1+1),\ldots,(j_r,j_r+1)\}
$$ 
be an arbitrary set of simple roots, where 
$1\leq j_1<\cdots<j_r\leq n$. For $0\leq i\leq r$, 
let $V_i$ be the subspace of $k^{n+1}$ spanned by
$v_{j_i+1},\ldots,v_{j_{i+1}}$ and let 
$V_{\leq i}=V_0\oplus\cdots\oplus V_i$, 
$V_{\geq i}=V_i\oplus\cdots\oplus V_r$. Then 
$$
k^{n+1}=V=V_0\oplus\cdots\oplus V_r,
$$ 
and $P_J$ (resp.~$P_J^-$ ; $L_J$) is the stabilizer in 
${\rm PGL}(V)$ of all subspaces $\mP(V_{\leq i})$
(resp. $\mP(V_{\geq i})$; $\mP(V_i)$) for $0\leq i\leq r$. One checks
that
$$
^J\!W\!^I=\{w_1,w_{j_1+1},\ldots,w_{j_r+1}\}
$$ 
and that $L_J\cap P^{w_i}$ is the stabilizer in $L_J$ of the
line $k v_{j_i+1}$, for every $w_i$ in $^J\!W\!^I$. 

Now the irreducible components of $p(\pi^{-1}(x_J))$ are the
$$
Z_i=\{(x,y)\in\mP(V)\times\mP(V) ~\vert~ x\in \mP(V_{<i}+\ell),
~y\in\mP(V_{>i}+\ell) \text{ for some line } \ell \text{ in } V_i\},
$$
for $0\leq i\leq r$, where $Z_i=Z(w_{j_i+1})$. Note that $Z_i$ meets
$Z_{i+1}$ along their divisor $\mP(V_{\leq i})\times\mP(V_{>i})$. 
We have, with obvious notation,
$$
\tilde Z_i=\{(x,y,\ell)\in\mP(V)\times\mP(V)\times\mP(V_i)~\vert~
x\in \mP(V_{<i}+\ell), ~y\in\mP(V_{>i}+\ell)\},
$$
and $\rho_i:\tilde Z_i\to Z_i$ is the first projection; note that
$\tilde Z_i$ is non-singular.

If $\dim(V_i)\geq 2$ and $0<i<r$, then the exceptional locus of
$\rho_i$ is 
$$
\mP(V_{<i})\times \mP(V_{>i})\times \mP(V_i).
$$ 
This subset has codimension $2$ in $\tilde Z_i$, and is mapped by
$\rho_i$ to $\mP(V_{<i})\times \mP(V_{> i})$. It follows that $Z_i$ is
singular along $\mP(V_{<i})\times \mP(V_{> i})$. On the other hand, if
$\dim(V_i)=1$, then 
$\tilde Z_i\cong Z_i=\mP(V_{\leq i})\times\mP(V_{\geq i})$
is non-singular. Finally, $Z_0$ (resp.~$Z_r$) is isomorphic to the
blow-up of $\mP(V_{>0})$ (resp.~of $\mP(V_{<r})$) in $\mP(V)$.

In particular, the ``minimal'' degenerations of the diagonal in
$\mP(V)\times\mP(V)$ correspond to the decompositions 
$V=V_0\oplus V_1$. Each such degeneration has $2$ irreducible
components: the blow-up $Z_0$ of $\mP(V_1)$ in $\mP(V)$, viewed as a
subvariety of $\mP(V_0)\times\mP(V)$, and the blow-up $Z_1$ of
$\mP(V_0)$ in $\mP(V)$, viewed in $\mP(V)\times\mP(V_1)$. These
components are glued along their common exceptional divisor
$\mP(V_0)\times\mP(V_1)$.

\medskip

\noindent
{\sl Remark.} If $X$ is the full flag variety, then every fiber of
$\pi$ is Gorenstein (for $\oB$ is Gorenstein, by \cite{BP} \S 5). But
this fails e.g. for $X=\mP^n$ where $n$ is even: let indeed $Z$ be the
total degeneration of the diagonal in $\mP^n\times\mP^n$. One checks
easily that the restriction
${\rm Pic}(\mP^n\times\mP^n)\to {\rm Pic}(Z)$
is an isomorphism. Assuming that $Z$ is Gorenstein, its canonical sheaf
is thus isomorphic to $\cO_Z(p,q)$ for unique integers $p$ and $q$. By
symmetry, we have $p=q$; by duality, we obtain
$$
\chi(Z,\cO_Z(-m,-m))= \chi(Z,\cO_Z(m+p,m+p))
$$
for all integers $m$. Moreover, the Hilbert polynomials of $Z$ and of
$\mP^n$ are equal, so that 
$$
\chi(Z,\cO_Z(m,m))=\chi(\mP^n,\cO_{\mP^n}(2m))=
\frac{(2m+1)\cdots(2m+n)}{n!}.
$$
Since $n$ is even, this yields a contradiction with the previous
equality.

\section{Wonderful completions via Chow varieties}

In this section, we assume that the characteristic of the ground
field $k$ is zero. As in Section 1, we consider a projective algebraic
variety $X$ together with a closed connected subgroup $G$ of 
$\Aut(X)$. We shall modify the constructions of Section 1 by
replacing Hilbert schemes with Chow varieties; the latter can be
defined as follows.

Choose a closed embedding of $X$ into some projective space $\mP$. Let
$\Chow_{n,d}(X,\mP)$ be the set of Chow forms of effective cycles of
$X$ with dimension $n$ and degree $d$. Then $\Chow_{n,d}(X,\mP)$ is a
projective algebraic set; moreover, the disjoint union of the 
$\Chow_{n,d}(X,\mP)$ over all $(n,d)$ is independent on the
projective embedding of $X$ (see \cite{Ba}, Corollaire, p.~115; this
may fail in positive characteristics, see \cite{Ko} I.4). 
We call this union the Chow variety of $X$ and denote it 
$\Chow(X)$; this definition differs from that in \cite{Ko} I.3, where
the seminormalization of $\Chow(X)$ is considered.

Taking the fundamental class of a subscheme defines the Hilbert-Chow
morphism, from the Hilbert scheme of $X$ (endowed with its reduced
scheme structure) onto the Chow variety; see \cite{Ko} I.3.15, I.3.23.3.

Let $\cC_{X,G}$ be the $G\times G$-orbit closure of the diagonal in
$\Chow(X\times X)$. This is a projective equivariant completion of
$G$; its points are graphs of elements of $G$, together with their
limits as cycles. The Hilbert-Chow morphism from $\Hilb(X\times X)$ to
$\Chow(X\times X)$ restricts to a $G\times G$-equivariant birational
morphism from $\cH_{X,G}$ to $\cC_{X,G}$. Moreover, the action of
$G\times G$ on $\cC_{X,G}$ can be interpreted in terms of the
(partially defined) composition of correspondences; for the latter,
see \cite{F} 16.1.

Consider now a semisimple adjoint group $G$ and a parabolic subgroup
$P$ such that $G$ acts faithfully on $G/P=X$. Then we saw that all
closed points of $\cH_{X,G}$ are reduced; thus, the Hilbert-Chow
morphism restricts to a bijective morphism 
$$
HC:\cH_{X,G}\to\cC_{X,G}.
$$ 
Composing with $\varphi:\oG\to\cH_{X,G}$, we obtain a morphism
$$
\psi:\oG\to\cC_{X,G},
$$
mapping every $x\in\oG$ to the cycle $p_*(\pi^*[x])$,
where all multiplicities equal $1$. We shall show that $\psi$ is an
isomorphism; by Theorem \ref{main}, this amounts to showing that $HC$
is an isomorphism. Actually, we shall obtain a direct proof of the
following equivalent statement.

\begin{theorem}\label{chow}
Let $X$ be a projective variety, homogeneous under a connected group
$G$ of automorphisms. Then $\cC_{X,G}$ is equivariantly isomorphic to
$\oG$.
\end{theorem}

\begin{proof} We begin with the following observation.

\begin{lemma}\label{closed} 
With the preceding notation and assumptions, every connected component
of $\Chow(X)$ contains a unique closed $G$-orbit, and admits a
$G$-equivariant embedding into the projectivization of a $G$-module.
\end{lemma}

\begin{proof} 
Let $\cC$ be a connected component of $\Chow(X)$. Then all points of
$\cC$, viewed as cycles in $X$, are algebraically equivalent. But
algebraic equivalence in $X$ coincides with rational equivalence, and
the Chow group of $X$ is freely generated by the classes of the
Schubert varieties (see \cite{F} Example 19.1.11). Therefore,
every cycle in $X$ is algebraically equivalent to a unique cycle with
$B$-stable support. Thus, $\cC$ contains a unique fixed point of $B$,
and hence a unique closed $G$-orbit.

For the second assertion, choose an equivariant embedding of $X$ into
the projectivization $\mP$ of a $G$-module. Then $\cC$ is contained in
some $\Chow_{n,d}(X,\mP)$. The latter is contained in the
projectivization of a $G$-module, by the construction of the
Chow variety.
\end{proof}

Applying Lemma \ref{closed} to $X\times X$, we see that $\cC_{X,G}$
contains a unique closed $G\times G$-orbit. The latter is isomorphic
to $G/B\times G/B$, by Lemma \ref{bijective}. Thus, Theorem \ref{chow}
is a consequence of

\begin{lemma}\label{miracle}
Let $V$ be a $G\times G$-module. Let $x\in\mP(V)$ satisfy the
following conditions:

\noindent
(i) The isotropy group $\Stab_{G\times G}(x)$ equals $\diag(G)$.

\noindent
(ii) The orbit closure $\overline{(G\times G)\cdot x}$ contains a
unique closed $G\times G$-orbit, and the latter is isomorphic to 
$G/B\times G/B$.

\noindent
Then the map 
$$
G\cong (G\times G)/\diag(G)\to\overline{(G\times G)\cdot x},~
(g,h)\mapsto (g,h)\cdot x
$$
extends to an isomorphism 
$\oG\to\overline{(G\times G)\cdot x}$.
\end{lemma}

\begin{proof} By \cite{Kn} \S 7, the map 
$G\to\overline{(G\times G)\cdot x}$ extends to an equivariant
birational morphism
$$
\oG\to\overline{(G\times G)\cdot x}.
$$
We shall construct an inverse to that morphism.

By (ii), there exists a unique line $\ell$ in $V$ such that the
corresponding point of $\mP(V)$ belongs to 
$\overline{(G\times G)\cdot x}$ and has isotropy group $B\times B$. 
Thus, $\ell$ consists of eigenvectors of $B\times B$ of weight
$(\lambda,\mu)$, where $\lambda$ and $\mu$ are regular dominant
weights of $B$. Let $V_{\lambda,\mu}$ be the $G\times G$-submodule
of $V$ generated by $\ell$. By complete reductibility, we may choose a 
$G\times G$-equivariant projection $V\to V_{\lambda,\mu}$. Then the
corresponding rational map 
$$
f:\mP(V)- - \to \mP(V_{\lambda,\mu})
$$ 
is $G\times G$-equivariant and defined at $\ell$, and hence defined
everywhere on $\overline{(G\times G)\cdot x}$. The image of $x$ under
$f$ is a fixed point of $\diag(G)$ in $\mP(V_{\lambda,\mu})$. In
particular, the $G\times G$-module $V_{\lambda,\mu}$ contains an
eigenvector of $\diag(G)$. Thus, this module is the space of
endomorphisms of the simple $G$-module $V_{\mu}$ with highest weight
$\mu$; moreover, the image of $x$ in 
$\mP(V_{\lambda,\mu})=\mP{\rm End}(V_{\mu})$ is the line spanned by
the identity map. By \cite{DP} 3.4, the $G\times G$-orbit closure of
that line is isomorphic to $\oG$. Thus, $f$ restricts to an
equivariant birational morphism from $\overline{(G\times G)\cdot x}$
onto $\oG$. 
\end{proof}

\end{proof}

Together with Zariski's main theorem, Theorem \ref{chow} implies that
the equivariant birational morphisms
$HC:\cH_{X,G}\to\cC_{X,G}$ and $\varphi:\oG\to\cH_{X,G}$ are
isomorphisms as well. This yields an alternative proof of Theorem
\ref{main} in characteristic zero; its only overlap with the proof of
Section 2 is the simple Lemma \ref{bijective}.

\medskip

Finally, we extend Theorem \ref{chow} to symmetric spaces. As at the
end of Section 2, consider an automorphism $\sigma$ of order $2$ of
$G$, the corresponding isomorphism $f:G/P\to G/\sigma(P)$, and the
closure $\cC_{f,G}$ of the $G$-orbit of the graph $\Gamma_f$ in the Chow
variety of $G/P\times G/\sigma(P)$. Then Lemma \ref{fixed} and Theorem
\ref{chow} imply readily

\begin{corollary}\label{chowtwist}
With preceding notation, $\cC_{f,G}$ is isomorphic to the wonderful
completion of $G/G^{\sigma}$.
\end{corollary}


\begin{thebibliography}{100}

\bibitem{A} D.~Akhiezer: Lie group actions in complex analysis,
Aspects of Math. E 27, Viehweg 1995.

\bibitem{Ba} D.~Barlet: Espace analytique r\'eduit des cycles
analytiques complexes compacts d'un espace analytique complexe de
dimension finie. In: Fonctions de plusieurs variables complexes II,
1-158, Lecture Note in Math. {\bf 482}, Springer-Verlag 1975.

\bibitem{Bo} N.~Bourbaki: Groupes et alg\`ebres de Lie, Chap.\,4,5,6,
Masson 1981.

\bibitem{B} M.~Brion: The behaviour at infinity of the Bruhat
decomposition, {\sl Comment. Math. Helv.} {\bf 73} (1998), 137-174.

\bibitem{BP} M.~Brion and P.~Polo: Large Schubert varieties,
{\sl Representation Theory} {\bf 4} (2000), 97-126.

\bibitem{De} M.~Demazure: Automorphismes et d\'eformations des
vari\'et\'es de Borel, {\sl Invent. math.} {\bf 39} (1977), 179-186. 

\bibitem{DP} C.~De Concini and C.~Procesi: Complete symmetric
varieties. In: Invariant Theory, 1-44, Lecture Note in Math. 
{\bf 996}, Springer-Verlag 1983.

\bibitem{DS} C.~De Concini and T.~A.~Springer: Compactification of
symmetric varieties, {\sl Transform. Groups} {\bf 4} (1999),
273-300. 

\bibitem{F} W.~Fulton: Intersection Theory, Ergeb. der Math.
{\bf 2}, Springer-Verlag 1998.

\bibitem{K1} M.~M.~Kapranov: Chow quotients of Grassmanians. I, 
{\sl Adv. Soviet Math.} {\bf 16} (1993), 29-110.

\bibitem{K2} M.~M.~Kapranov: Veronese curves and Grothendieck-Knudsen
moduli space $\overline{M}_{0,n}$, {\sl J. Alg. Geom.} {\bf 2} (1993),
239-262.

\bibitem{Kn} F.~Knop: The Luna-Vust theory of spherical embeddings,
in: Proceedings of the Hyderabad Conference on Algebraic
Groups, 225-250, Manoj Prakashan 1991.

\bibitem{Ko} J.~Koll\'ar: Rational curves on algebraic varieties,
Ergeb. Math. {\bf 32}, Springer-Verlag 1996.

\bibitem{L} L.~Lafforgue: Pavages des simplexes, sch\'emas de graphes
recoll\'es et compactification des ${\rm PGL}^{n+1}_r/{\rm PGL}_r$,
{\sl Invent. math.} {\bf 136} (1999), 233-271.

\bibitem{LP} P.~Littelmann and C.~Procesi: Equivariant cohomology of
wonderful compactifications, in: Operator algebras, unitary
representations, enveloping algebras, and invariant theory (Paris,
1989), 219-262, Progr. Math. {\bf 92}, Birkh\"auser 1990.

\bibitem{Sp} T.~A.~Springer: Linear algebraic groups, Progress in
Math. {\bf 9}, Birkh\"auser 1998.

\bibitem{Sn} D.~Snow: Transformation groups of compact K\"ahler
spaces, {\sl Arch. Math. (Basel)} {\bf 37} (1981), 364-371.

\bibitem{St} E.~Strickland: A vanishing theorem for group
compactifications, {\sl Math. Ann.} {\bf 277} (1987), 165-171.

\bibitem{T} M.~Thaddeus: Complete collineations revisited, 
{\sl Math. Ann.} {\bf 315} (1999), 489-495.

\end{thebibliography}
\end{document}